\documentclass[11pt]{amsart}
\numberwithin{equation}{section}
\newtheorem{theorem}{Theorem}

\newtheorem{proposition}{Proposition}
\newtheorem{lemma}{Lemma}
\newtheorem{corollary}{Corollary}
\numberwithin{theorem}{section} \numberwithin{lemma}{section}
 \numberwithin{definition}{section}
\numberwithin{proposition}{section}
\newtheorem{remark}{Remark}[section]
\def\R{\bf R}

\def\al{\aligned}
\def\eal{\endaligned}
\def\M{{\bf M}}
\def\be{\begin{equation}}
\def\ee{\end{equation}}
\def\lab{\label}

\def\R{\bf R}
\def\M{\bf M}

\def\al{\aligned}
\numberwithin{equation}{section}

\begin{document}

\tracingpages 1
\title[backward-control]{\bf A formula for backward and control problems of the heat equation }
\author{Qi S. Zhang}
\address{ Department of
Mathematics, University of California, Riverside, CA 92521, USA }
\date{May 2020: key words: heat equation, exact control, explicit formula. MSC2020: 35k05, 49K20, 58J35
}

\begin{abstract}

(a). Using time analyticity result, we address a basic question  for a nonhomogeneous backward heat equation (exact control problem) in the setting of smooth domains and compact manifolds, namely: when is essentially time independent control  possible? i.e. The control function is 0 on one time interval and stationary on the other. For general $L^2$ initial values, the answer is: if and only if the full space domain is used for the control function.  Also an explicit formula for the control function is found in the form of an infinite series involving the heat kernel, which converges rapidly.

(b). A formal exact formula for a time dependent control function supported in a proper subdomain is also obtained via eigenfunctions of the Laplacian. The formula is rigorous on any finite dimensional space spanned by the eigenfunctions and there is no smoothness assumption on the whole domain, making partial progress on a problem on p74 \cite{Zu:1}.

(c). A byproduct is an inversion formula for the heat kernel.

\end{abstract}
\maketitle
\section{Introduction}

In this paper we consider a backward nonhomogeneous problem of the heat equation, which belongs to the following typical control problem involving an evolution equation. Given an initial state or value, can one find a nonhomogeneous or control term to reach
a desired final state or value in given time?

Here is a sample of the main results in the classical papers \cite{LR:1} (Corollary 1) or \cite{FI:1}, after some translation and recast.

{\it Let $D$ be a smooth, bounded domain in a Riemannian manifold $\M$, $G=G(x, t, y)$ be the heat kernel on $D$ with zero boundary value. Given any $u_0, u_1\in L^2(D)$ and a subdomain $\omega \subset D$, suppose
\be
\lab{zgu1}
z=z(x)=\int_D G(x, t, y) u_1(y) dy,
\ee  then there exists a  function $g \in L^2(D \times (0, T])$ such that
the problem below has a classical solution.
\be
 \lab{cp0}
 \begin{cases}
 \Delta u(x, t) - \partial_t u(x, t)=g(x, t) \chi_\omega, \quad (x, t) \in D \times (0, T], \\
 u(x, t) = 0, \quad (x, t) \in \partial D \times (0, T],\\
 u(x, T)= z(x), \\
 u(x, 0) = u_0(x).
 \end{cases}
 \ee
}  In other words, if a state can be reached by the free heat flow, then it can be reached by control from any $L^2$ initial state. In particular, $0$ state can always be reached i.e. exact null control is always possible.  This result has been widely extended and generalized and it has stimulated much further research. See \cite{LZ:1}, \cite{Zu:1}, \cite{LL:1},  \cite{LT:1}, \cite{LZZ:1}, \cite{Tr:1}, \cite{DM:1}, \cite{EMZ:1}, \cite{ABGM:1} and the references therein e.g. See also related earlier work \cite{Lio:1},  \cite{FR:1}, \cite{Ru:1},
 \cite{Lp:1} and \cite{RW:1}.
Nevertheless, there are some basic questions which are important from both theoretical and practical point of view.

Question 1:   How does one characterize all states that can be reached by the free heat flow? i.e.
When does (\ref{zgu1}) hold for $z=z(x)$?

Question 2: How does one determine the control function $g$?

Question 1, being equivalent to the solvability problem of the backward heat equation, has been answered explicitly in a recent paper \cite{DZ:1} in full generality.  See also \cite{Z:1} under a worse condition and a generalization in \cite{DP:1}. In case of bounded domains with zero boundary values, a necessary and sufficient condition involving all eigenfunctions of $\Delta$ is found in \cite{Lin:2} p249 and  see also an abstract criteria \cite{CJ:1} Theorem 9.

\begin{corollary} (\cite{DZ:1})
\lab{cobhe}
Let $\M$ be a complete, $d$ dimensional, noncompact Riemannian
 manifold such that the Ricci curvature satisfies $Ric \ge
 -(d-1) K_0$ for a nonnegative constant $K_0$. Then the Cauchy problem for the
backward heat equation
\be
\lab{bhe}
\begin{cases} \partial_t u+\Delta u  = 0,\\
u(x, 0)=z(x)
\end{cases}
\ee has a smooth  solution of exponential growth of order 2 in time interval
$(0,
\delta)$ for some $\delta>0$ if and only if
\be
\lab{aijie}
|\Delta^j z(x)| \le A^{j+1}_3 j^j e^{ A_4  d^2(x, 0)}, \quad
j=0, 1,2,\ldots,
\ee where  $A_3$ and $A_4$ are some positive constants.
\end{corollary}

In the above, the meaning that $u=u(x, t)$ is of exponential growth of order $2$ is that
$|u(x, t)| \le C_1 e^{C_2 d^2(x, 0)}$ for all $x \in \M$ and $t$ in some given interval. Here $0 \in \M$
is a reference point and $C_1, C_2$ are positive constants. This condition is sharp due to Tychonov's non-uniqueness example which can be extended to the backward heat equation by reflection in time.

Although the above result is stated for a noncompact manifold, as mentioned in that paper, the conclusion still holds for compact manifolds and for smooth domains with $0$ boundary condition, and the proof is simpler. Moreover the exponential term and the curvature condition all become redundant. So in the compact setting a state $z=z(x)$ is reachable by the free heat flow if and only if
\be
\lab{areach}
|\Delta^j z(x)| \le C C^j_* j!, \quad j=0, 1, 2, ...
\ee for some positive constants $C, C_*$. Here $C_*$ may depend on the length of the time interval. See Lemma \ref{lebhe} below for an explicit estimate of $C_*=e^+/T$ where $T$ is the length of time.

The goal of the current paper is to address a part of Question 2.

The classical variational method in \cite{FI:1} and the method in \cite{LR:1} provides an implicit way of finding a time dependent control function which may be supported in a given subdomain in space time. A minimizer function at a future time is obtained by minimizing certain functionals in a Hilbert space whose norm is defined in space time. Then the control function is determined as a cut off of the solution  of the backward heat equation with the minimizer as the final value.  It is hard to find the exact minimal value of the functionals involved and
 the control function is usually complicated and not smooth. See \cite{MZ:1} and \cite{EZ:1} e.g.

In applications, smooth and time independent controls are useful if they can be found.
 It turns out that, if the full space domain is used, then to reach the same final state as the time dependent controls, the control function (nonhomogeneous term) can be essentially independent of time, i.e. it is 0 on one time interval and stationary on the other. Moreover an explicit formula for the control function is found in the form of an infinite series involving the heat kernel. We also show that if the control function is supported in a proper subdomain, then this essentially time independent control is impossible in general.   What causes the difference?  Comparing with the traditional method of using weighted energy estimates (Carleman estimates), the new input is the time analyticity of solutions of the heat equation with stationary nonhomogeneous terms. This allows us to use power series in time to represent solutions and carry out calculations. See Section 2.

A formal exact formula for a control function supported in a proper subdomain is obtained via eigenfunctions of the Laplacian. The formula is rigorous and calculable on any finite dimensional spaces spanned by the eigenfunctions and there is no smoothness assumption on the whole domain. A remark is also made on the impossibility of null control with rough coefficients.  See Section 3 for details.

An  inversion formula for the heat kernel, as a by product,  is given in Section 4, which may be of independent interest for inverse problems, among others.

\section{A formula for a control function in the full domain}

In order to state the result, let us first introduce a bit of notations.
We use $\M$ to denote a $n$ dimensional,  Riemannian manifold, $\Delta$ is the Laplace-Beltrami operator, $G=G(x, t, y)$ to denote the heat kernel with $0$ boundary condition on domain $D$;  and $0$ a reference point on $\M$, $d(x, y)$ is the geodesic distance for $x, y \in \M$. We use $C, C_1, ...$ etc to denote positive constants, which may change in value; $e=\lim_{n \to \infty}(1+1/n)^n$. The manifold setting for the main result is chosen for convenience. One could also choose to work on certain metric spaces or for simplicity smooth bounded domains in $\R^n$.

\begin{theorem}
\lab{thmain}
 Let $\M$ be a n dimensional, compact Riemannian manifold without boundary and $D \subset \M$ be a smooth domain or $D=\M$.

 (a).
  Suppose any initial value $u_0 \in L^2(D)$ and a function $z \in C^\infty(D) \cap C_0(D)$ are given, which satisfies, for some positive constants $C, A$,
 \be
 \lab{djau00}
 |\Delta^j z(x) |  \le C A^j \,  j !, \qquad \forall x \in {\M}, j=0, 1, 2, ...
 \ee
 Let $\delta =\min \{ \frac{1}{2 A}, \frac{1}{1+ 2 e} T \}$. Then for any $T_0 \in (T-\delta, T]$, there exists a control function $f \in C^\infty(D) \cap C_0(D)$ such that the nonhomogeneous problem (control problem)
 \be
 \lab{tcp}
 \begin{cases}
 \Delta u(x, t) - \partial_t u(x, t)=f(x) \,  \chi_{[T_0, T]}(t), \quad (x, t) \in D \times (0, T], \\
 u(x, t) = 0, \quad (x, t) \in \partial D \times (0, T],\\
 u(x, T)= z(x), \\
 u(x, 0) = u_0(x),
 \end{cases}
 \ee has a continuous solution $u$ which is smooth except at $T_0$. Moreover $f$ is given by the formula
 involving the Dirichlet heat kernel $G$:
 \be
 \lab{tcf}
 f(x) =  \int_{D} \sum^\infty_{k=1} \Delta_y G(x, 2^k (T-T_0), y) \, (b(y)-u(y, T_0)) dy,
 \ee
 \be
 \lab{tfb}
 b=b(x)= \sum^\infty_{j=0} \Delta^j z(x) \frac{(T_0-T)^j}{j !}.
 \ee

In particular, if $z=z(x)$ is reachable by the free heat flow from initial time $0$ to $T$, i.e. (\ref{zgu1}) holds,  then the above conclusions hold with $\delta =\frac{T}{1+2 e}$.

(b). Let $D_0 \subset \overline{D_0} \subset D$ be a proper subdomain.  Then for any
$f \in L^1(D)$ supported in $D_0$, any initial value $u_0 \in C_0(D)$ which does not vanish in $D/D_0$,  the null control problem, i.e. problem (\ref{tcp}) with $u(\cdot, T)=0$ and $T_0=0$,
has no solution $u \in C((0, T], L^2(D))$.
\end{theorem}

{\remark When $D$ is the whole manifold $\M$, i.e. $\partial D = \emptyset$, then no boundary condition is imposed. The smoothness of $\partial D$ can be relaxed to $C^2$ condition.}

{\remark The idea of using free heat flow, i.e. no control in the first part of the time interval, seems natural due to the smoothing effect of the free heat flow. In the second part of the time interval, the
stationary control is easier to execute than time dependent control in practice. As shown in the proof, the series in (\ref{tcf}) converges rapidly, allowing practical computation. Similar results can be extended to some other evolution equations with time analyticity property.
}

{\remark
According to Lemma \ref{lebhe} below, if $z=z(x)$ is reachable by the free heat flow from $0$ to $T$, then (\ref{djau00}) holds with $A=e^+/T$ where $e^+$ is any number greater than $e=2.71828...$. So the theorem allows final states which are broader than those reachable by the free heat flow from $0$ to $T$ since there is no restriction to the size of $A$.
}

{\remark
Results of the theorem can be extended to the case when solutions $u$ and functions $u_0, z, f$ are Schwartz functions on some noncompact manifolds including $\R^n$.
}

Since the proof for the case $\partial D = \emptyset$ is almost identical to the case when $\partial D \neq \emptyset$, we will just concentrate on the former and indicate a few necessary changes in the proof.

We will need two lemmas before finishing the proof of Theorem \ref{thmain} at the end of the section. The first one is the main technical result of the paper. Here we solve the control problem when time is sufficiently short and the initial value and final value are in the same class.

\begin{lemma} (main lemma)
\lab{lemain}
 Let $\M$ be a n dimensional, compact Riemannian manifold without boundary and $D \subset \M$ a smooth domain.
 Suppose $z=z(x), \, u_0=u_0(x)$ are given smooth functions in $C^\infty(D) \cap C_0(D)$ such that
 \be
 \lab{djau0}
 |\Delta^j z(x) | + |\Delta^j u_0(x) | \le C C^j_* \,  j !
 \ee for all $x \in D$, $j=0, 1, 2, ...$. Here $\Delta$ is the Laplace-Beltrami operator, $C$ and $C_*$
 are positive constants. Then for any
 \be
 \lab{tjie} T \in (0, \frac{1}{2 C_*}),
 \ee there exists a control function $f \in C^\infty(D) \cap C_0(D)$ such that the nonhomogeneous problem (control problem)
 \be
 \lab{cp}
 \begin{cases}
 \Delta u(x, t) - \partial_t u(x, t)=f(x), \quad (x, t) \in D \times (0, T], \\
 u(x, t)=0, \quad (x, t) \in \partial D \times (0, T], \\
 u(x, T)= z(x), \\
 u(x, 0) = u_0(x),
 \end{cases}
 \ee has a unique smooth solution $u$. Moreover $f$ is given by the formula
 \be
 \lab{cf}
 f(x) =  \int_{D} \sum^\infty_{k=1} \Delta_y G(x, 2^k T, y) \, (b-u_0)(y) dy,
 \ee
 \be
 \lab{fb}
 b= \sum^\infty_{j=0} \Delta^j z(x) \frac{(-T)^j}{j !}.
 \ee In addition,
 \be
\lab{dif0}
| \Delta^i f(x)| \le C_5 \frac{C^{i+1}_* (i+1) !}{(1-C_* T)^{i+1}} (1+ \frac{1}{ T})
\ee with $C_5$ depends only on $\M$; and
  $u$ is given by the formula
 \be
 \lab{usol0}
u(x, t)= z(x) + \sum^\infty_{j=1} (\Delta^j z - \Delta^{j-1} f)(x) \frac{ (t-T)^j}{j !}.
\ee
 \proof
\end{lemma}

As mentioned, we will only give a proof for the case when $D = \M$, i.e. $\partial D = \emptyset$. Otherwise the proof is almost identical. One just needs to replace the term $G(x, t, y) - \frac{1}{|\M|}$ below by $G(x, t, y)$ and also to make sure boundary terms vanish in integration by parts.

The proof is carried out in 3 steps.

{\it Step  1.} We show that the functions $f$ in (\ref{cf}) and $b$ in (\ref{fb}) are well defined in the sense that the series converge absolutely and uniformly.

From the conditions (\ref{djau0}) and $0 \le T C_*<1$,
\[
|\Delta^j z(x) \frac{(-T)^j}{j !}| \le C ( C_* T)^j.
\]Hence the series in (\ref{fb}) converges absolutely and uniformly and  the function $b$ is well defined.
Similarly
\be
\lab{deltab}
\Delta b= \sum^\infty_{j=0} \Delta^{j+1} z(x) \frac{(-T)^j}{j !}
\ee which also converges absolutely and uniformly.

Next we prove that $f$ is well defined by heat kernel estimates. Since $\M$ is a compact manifold, it is well known that the following upper bound holds for the heat kernel $G=G(x, t, y)$: for all $x, y \in \M$ and $t>0$,
\be
\lab{hkup}
\left| G(x, t, y) - \frac{1}{|\M|} \right| \le C_1 \frac{e^{-C_2 t}}{t^{n/2}} e^{- C_3 d^2(x, y)/t}.
\ee Here $|\M|$ is the volume of the manifold $\M$,  $C_1, C_2, C_3$ are positive constants and $d(x, y)$
is the distance between $x$ and $y$. This bound can be found implicitly or explicitly in many references such as \cite{Li:1} Chapter 16 e.g..
Note that for large $t$ the right hand side of the bound decays exponentially due to the term $e^{- C_2 t}$. We mention that if $\M$ is noncompact, then there would be a generic, exponentially growing upper bound for $G$ in general. Using a mean value inequality and the property that $\partial_t G(x, t, y)$ is a solution of the heat equation for the variables $(x, t)$, it is not hard to deduce the bound
\be
\lab{hkupt}
\left| \partial_t G(x, t, y) \right| = \left| \Delta_x G(x, t, y) \right|=\left| \Delta_y G(x, t, y) \right| \le C_1 \frac{e^{-C_2 t}}{t^{(n+2)/2}} e^{- C_3 d^2(x, y)/t},
\ee for all $x, y \in \M$ and $t>0$. Here the positive constants $C_1, C_2, C_3$ may have changed. We will give a proof of this bound below,
since, comparing with generic bounds,  we need the large time decay property.

Let us start with a well known parabolic mean value inequality which can be found in Theorem 14.7 of \cite{Li:1} e.g.  Suppose $u$ is a positive subsolution to the heat equation on ${\M} \times [0, T]$. Let $T_1, T_2 \in [0, T]$ with $T_1<T_2$, $R>0$, $p>0$, $\delta, \eta \in (0, 1)$. Then there exist  positive constants $C_1$ and $C_2$, depending only on $p, n$ such that
\be
\lab{pmvi}
\al
\sup_{B(0, (1-\delta) R) \times [T_1, T_2]} u^p &\le C_1 \frac{\bar{V}(2R)}{|B(0, R)|}
( R \sqrt{K_0} + 1) \exp ( C_2 \sqrt{K_0 (T_2-T_1)}) \\
&\qquad \times
\left(\frac{1}{\delta R} + \frac{1}{\sqrt{\eta T_1}} \right)^{n+2} \,
\int^{T_2}_{(1-\eta) T_1} \int_{B(0, R)} u^p(y, s) dy ds;
\eal
\ee Here $\bar{V}(R)$ is the volume of geodesic balls of radius $R$ in the simply connected space form with constant sectional curvature $-K_0$; $|B(0, R)|$ is the volume of the geodesic ball $B(0, R)$ with center $0$ and radius $R$.  This mean value inequality is stated for a general complete manifold. A small difference needs to be mentioned in that the term $\exp ( C_2 \sqrt{K_0 (T_2-T_1)})$ here is $\exp ( C_2 \sqrt{K_0 T_2})$ in \cite{Li:1}. But the proof is identical after doing a time shift $T \to T-T_1$.

  In our case
$\M$ is compact, $R$ can be chosen as one half of the diameter of $\M$, $T_2 =T$ and $T_1 = \sup \{ T/2, T-1\}/(1-\eta)$.
With these choices, writing $u=u(x, t)= G(x, t, y)$, we infer from (\ref{pmvi}) that
\be
\lab{mviq}
\sup_{{\M} \times [T_3, \, T]} (\partial_t u)^2 \le \frac{C_4}{\min \{T^{(n+2)/2}, 1 \}}
\int^{T}_{(1-\eta) T_1} \int_{\M} (\partial_s u)^2(y, s) dy ds;
\ee with $T_3=\sup\{2T/3, \,  T-1/2\}$. Here $C_4$ depends on $\M$ through $K_0$, $|\M |$ and the diameter of $\M$.

  Denote by $\psi=\psi(t)$ a standard Lipschitz cut off function supported in $[(3T_3 + T)/4, T]$ such that $\psi=1$ in $[(T_3+T)/2, T]$ and $|\partial_t \psi| \le C/\min\{T, 1\}$. Since $u$ is a smooth solution to the heat equation, we can compute, denoting $Q= {\M} \times [T_3, T]$ and $\partial_t u =u_t$, that
\[
\al
\int_{Q} &(\Delta u)^2 \psi^2 dxdt = \int_{Q} u_t \Delta u \psi^2 dxdt\\
&=-\int_{Q} ((\nabla u)_t \nabla u) \, \psi^2 dxdt
= - \frac{1}{2} \int_{Q} (|\nabla u |^2)_t \, \psi^2 dxdt \\
&\le \frac{1}{2} \int_{Q} |\nabla u |^2 \, (\psi^2)_t dxdt.
\eal
\]Therefore
\[
\int_{Q} (\Delta u)^2 \psi^2 dxdt \le \frac{C}{(\min\{T, 1\})} \int^T_{(3 T_3+T)/4} \int_{\M} |\nabla u |^2 dxdt.
\]This and the standard Cacciopoli inequality (energy estimate)  show that
\be
\lab{ddu2r-4}
\int_{Q} (\partial_t u)^2 \psi^2 dxdt =\int_{Q} (\Delta u)^2 \psi^2 dxdt \le \frac{C}{(\min\{T, 1\})^2} \int_{Q} \left|u-\frac{1}{|\M|} \right|^2  dxdt.
\ee Now (\ref{hkupt}) follow from (\ref{ddu2r-4}), (\ref{mviq}) and (\ref{hkup}). An immediate consequence of (\ref{hkupt}) is
\be
\lab{g2jt}
|\Delta_y G(x, 2^j T, y)| \le \frac{C_1  \exp ( - C_2 2^j T)}{( 2^j T)^{(n+2)/2}}.
\ee Therefore the series
\[
\sum^\infty_{k=1} \Delta_y G(x, 2^k T, y)
\]converges uniformly and absolutely for each fixed $T>0$ since it is dominated by
\[
\sum^\infty_{k=1} \frac{C_1  \exp ( - C_2 2^j T)}{( 2^j T)^{(n+2)/2}}.
\]This implies that the function $f$ in (\ref{cf}) is well defined.
\medskip

{\it Step 2}.  In this step, we derive bounds on the functions $|\Delta^i f|$ (\ref{dif}) and $|\Delta^i b|$ (\ref{dibb}), $i=0, 1, 2, ...$.

For any positive integer $i$, by Assumption (\ref{djau0}) and Definition (\ref{fb}), we have
\be
\lab{dib}
\al
|\Delta^i b| &= \left |  \sum^\infty_{j=0} \Delta^{i+j} z(x) \frac{(-T)^j}{j !} \right |\\
&\le C C^i_* \sum^\infty_{j=0} (i+j) ...(1+j) ( C_* T)^j\\
&= C C^i_* \partial^i_s \sum^\infty_{j=0} s^j, \qquad \text{with} \quad s \equiv C_* T\\
&= C C^i_* \partial^i_s (1/(1-s)).
\eal
\ee This shows
\be
\lab{dibb}
|\Delta^i b| \le C \frac{ C^i_*  i!}{(1 - C_* T)^{i+1}}.
\ee

Next we derive a bound for $|\Delta^i f|$. From (\ref{cf}), using $\Delta_x G(x, t, y)= \Delta_y G(x, t, y)$, we compute
\be
\lab{dicf}
\al
 \Delta^i f(x) &=  \int_{\M} \sum^\infty_{k=1} \Delta^i_x \Delta_y G(x, 2^k T, y) \, (u_0-b)(y) dy\\
&=  \int_{\M} \sum^\infty_{k=1}  \Delta^{i+1}_y G(x, 2^k T, y) \, (u_0-b)(y) dy \\
&= \int_{\M} \sum^\infty_{k=1}  \left( G(x, 2^k T, y) - \frac{1}{|\M|}
 \right)\, \Delta^{i+1} (u_0-b)(y) dy,
\eal
\ee where integration by parts was performed in the last step. Therefore
\[
| \Delta^i f(x)| \le   \sup_y |\Delta^{i+1} (u_0-b)(y)| \, \int_{\M} \sum^\infty_{k=1}  \left| G(x, 2^k T, y) - \frac{1}{|\M|}
 \right|  dy,
\] which infers, by the assumed bounds on $u_0$ (\ref{djau0}), bounds on $b$ (\ref{dibb}) and the heat kernel bound (\ref{hkup}) that
\[
| \Delta^i f(x)| \le  C_1 \frac{C^{i+1}_* (i+1) !}{(1-C_* T)^{i+1}}    \sum^\infty_{k=1} \frac{ e^{-C_2 2^k T}}{( 2^k T)^{n/2}}
\int e^{- C_3 d^2(x, y)/(2^k T)} dy.
\] By direct computation
\[
\frac{1}{( 2^k T)^{n/2}}
\int e^{- C_3 d^2(x, y)/(2^k T)} dy \le C_4
\]where $C_4$ is a positive constant depending on the on the volume lower bound of $\M$ and its volume doubling constant.  The last two inequalities together imply the bound
\be
\lab{dif}
| \Delta^i f(x)| \le 2 C \frac{C^{i+1}_* (i+1) !}{(1-C_* T)^{i+1}} (1+ \frac{C_1 C_4}{C_2 T}).
\ee Note that $C_1, C_2, C_4$ depend only on $\M$. Thus (\ref{dif0}) holds.

\medskip

{\it Step 3.} Now we prove that (\ref{cp}) has a solution given by

\be
\lab{useri}
u(x, t)= \sum^\infty_{j=0} a_j(x) \frac{ (t-T)^j}{j !},
\ee with $a_0(x)=z(x)$ the final value and $a_j$, $j=1, 2, 3, ...$ given by the recurrence formula:
\be
\lab{huigui}
\begin{cases}
\Delta a_0 = a_1 + f, \\
\Delta a_j = a_{j+1}, \qquad j=1, 2, 3, ...
\end{cases}
\ee i.e.
\[
a_1 = \Delta a_0 - f, \quad ..., \quad a_j = \Delta^j a_0 - \Delta^{j-1} f, ....
\]

Differentiating (\ref{useri}) formally, one obtains
\be
\lab{seridu}
\Delta u =  \sum^\infty_{j=0} \Delta a_j(x) \frac{ (t-T)^j}{j !},
\ee and
\be
\lab{seridtu}
\partial_t u =  \sum^\infty_{j=0}  a_{j+1}(x) \frac{ (t-T)^j}{j !}.
\ee In order that $u$ satisfies the nonhomogeneous heat equation in (\ref{cp}), the recurrence formula (\ref{huigui}) must hold.
To complete the proof, we need to do two things. One is to prove the three series (\ref{useri}),
(\ref{seridu}) and (\ref{seridtu}) all converge absolutely and uniformly. Two is to show that $u(x, 0)=u_0(x)$.

We will just prove that the series
\be
\lab{usol}
u(x, t)= a_0(x) + \sum^\infty_{j=1} (\Delta^j a_0 - \Delta^{j-1} f) \frac{ (t-T)^j}{j !}
\ee converges uniformly and absolutely. The other two series can be handled similarly.
By the bounds (\ref{dif}) for $\Delta^{j-1} f$ and the assumed bound (\ref{djau0}) on $\Delta^j a_0$,
we see that the general term in the series (\ref{usol}) is dominated by:
\[
| \Delta^j a_0 - \Delta^{j-1} f| \frac{ (t-T)^j}{j !}
\le C_5 (1+ \frac{1}{T}) \frac{C^j_* j! T^j}{(1-C_* T)^j j!} =  C_5 (1+ \frac{1}{T}) \left(\frac{C_* T}{1-C_* T} \right)^j.
\] The assumption that $T \in (0, \frac{1}{ 2 C_*})$ gives
\[
\frac{C_* T}{1-C_* T}<1.
\]Hence the series (\ref{usol}) converges uniformly and absolutely.

So we are left to verify $u(x, 0)=u_0(x)$.  This is done by setting $u(x, 0)=u_0(x)$ in (\ref{useri}) and prove that it is equivalent to $f$ being given by (\ref{cf}). Taking $t=0$ in (\ref{useri}) and suppose $u(x, 0)=u_0(x)$ as desired.
Then
\[
u_0(x) = a_0(x) + \sum^\infty_{j=1} (\Delta^j a_0 - \Delta^{j-1} f)(x) \frac{ (-T)^j}{j !},
\]which yields
\be
\lab{fbu0}
\al
\sum^\infty_{j=1}  \Delta^{j-1} f(x) \frac{ (-T)^j}{j !}& = a_0(x) + \sum^\infty_{j=1} \Delta^j a_0(x) \frac{ (-T)^j}{j !} -u_0(x) \\
&= b(x) - u_0(x).
\eal
\ee Here we have used the definition of $b(x)$ in (\ref{fb}). Note the series converges uniformly and absolutely as proven in the previous paragraph.

Next we will  invert (\ref{fbu0}) so that $f$ will appear explicitly.  This is the key argument in this proof.
Multiplying the heat kernel $G=G(x, 2T-t, y)$ on both sides of (\ref{fbu0}) and integrate, we find that
\be
\lab{intgdjf}
\sum^\infty_{j=1} \int G(x, 2 T-t, y) \Delta^{j-1} f(y) dy \frac{ (-T)^j}{j !} = \int G(x, 2 T- t, y) ( b-u_0)(y) dy.
\ee Using
\[
\Delta_y G(x, 2T-t, y) + \partial_t G(x, 2T-t, y) = 0, \quad t< 2T,
\] and integration by parts, we see, since $\Delta^j f = \Delta^j G=0$ on the boundary whenever it is nonempty, that
\[
\al
\int G(x, 2T - t, y) \Delta^{j-1} f(y) dy & = \int \Delta^{j-1}_y G(x, 2 T - t, y)  f(y) dy\\
&= (-1)^{j-1} \partial^{j-1}_t \int G(x, 2 T - t, y)  f(y) dy.
\eal
\] Substituting this to (\ref{intgdjf}) we arrive at the identity
\be
\lab{intgdjf2}
\sum^\infty_{j=1} \partial^{j-1}_t \int G(x,2 T- t, y)  f(y) dy \,  \frac{ T^j}{j !} =  - \int G(x, 2 T- t, y) ( b-u_0)(y) dy.
\ee
Writing
\be
\lab{ABxt}
A(x, t) =  \int G(x, 2 T - t, y)  f(y) dy, \qquad B(x, t)= - \int G(x, 2 T - t, y) ( b-u_0)(y) dy,
\ee then (\ref{intgdjf2}) can be written as
\be
\lab{seriAB}
\sum^\infty_{j=1}   \frac{ T^j}{j !} \, \partial^{j-1}_t A(x, t) = B(x, t).
\ee

Observe that, for the variable $s= 2 T - t$ and $x$,  both $A(x, t)$ and $B(x, t)$ are solutions of the heat equation with bounded initial values.
According to Theorem 2.1 in \cite{DZ:1}, they are real analytic in time for all $s = 2 T -t>0$, i.e. $t < 2 T$.  We should mention that
that theorem was stated for noncompact manifolds with Ricci curvature bounded from below. However the conclusion is still valid for the current compact setting since the proof is actually simpler and without the need of spatial cut-off functions. Differentiating (\ref{seriAB}) with respect to $t$ gives
\be
\lab{tjab1}
\sum^\infty_{j=1}   \frac{ T^j}{j !} \, \partial^{j}_t A(x, t) =  \partial_t B(x, t).
\ee The convergence of the above series is justified due to the bounds (\ref{dif}) for $\Delta^{j-1} f$ since
\[
\al
\partial^{j}_t A(x, t)&= \partial^{j}_t \int G(x, 2 T - t, y)  f(y) dy = (-1)^j \int \Delta^j_y G(x, 2 T -t, y)  f(y) dy\\
& = (-1)^j \int  G(x, 2 T- t, y) \Delta^j_y f(y) dy
\eal
\] so that
\[
|\partial^{j}_t A(x, t)| \le 2 C \frac{C^{j+1}_* (j+1) !}{(1-C_* T)^{j+1}} (1+ \frac{C_1 C_4}{C_2 T}),  \quad
t \in (0, 2T).
\] Using this and  $C_* T/(1-C_* T)<1$ again, we know that (\ref{tjab1}) converges uniformly and absolutely for $t \in (0, 2 T)$.

Therefore (\ref{tjab1}) yields:
\[
\sum^\infty_{j=0}   \frac{ T^j}{j !} \, \partial^{j}_t A(x, t) - A(x, t) =  \partial_t B(x, t).
\]Since $A=A(x, t)$ is analytic in time for all $t<2 T$, Taylor expansion around $t$ with a fixed $t<T$ reads
\[
\sum^\infty_{j=0}   \frac{ T^j}{j !} \, \partial^{j}_t A(x, t) = A(t+T).
\]The two preceding identities imply the relation
\be
\lab{AtTB}
A(x, t+T)-A(x, t) = \partial_t B(x, t), \quad i.e.
\ee
\[
\al
\int & G(x, T -t, y) f(y) dy - \int G(x, 2 T -t, y) f(y) dy \\
&= - \partial_t \int G(x, 2 T-t, y) (b-u_0)(y) dy\\
&=  \Delta_x \int G(x, 2 T - t, y) (b-u_0)(y) dy =  \int \Delta_y G(x,  2 T - t, y) (b-u_0)(y) dy\\
&=  \int  G(x,  2 T - t, y) \Delta (b-u_0)(y) dy.
\eal
\] Taking $t \to T^-$,  since the heat kernel $G=G(x, T-t, y)$ converges to the Delta function, we deduce
\[
f(x) = \int  G(x, T, y) \Delta (b-u_0) (y) dy + \int  G(x, T, y) f(y) dy.
\] Iterating once, we reach
\[
\al
f(x) &= \int  G(x, T, y) \Delta (b-u_0) (y) dy \\
 &\quad + \int G(x, T, w) \int  G(w, T, y) \Delta (b-u_0) (y) dy dw
+ \int  G(x, T, w) \int G(w, T, y) f(y) dy dw.
\eal
\]By the reproducing formula for the heat kernel, this becomes
\[
\al
f(x) &= \int  G(x, T, y) \Delta (b-u_0) (y) dy \\
 &\quad + \int G(x, 2 T, y)  \Delta (b-u_0) (y) dy dw
+ \int  G(x, 2 T, y) f(y) dy.
\eal
\]
Repeating this process, we arrive at (\ref{cf}):
\[
\al
f(x) &= \int \sum^\infty_{k=1} G(x,  2^k T, y) \Delta (b-u_0) (y) dy\\
&= \int \sum^\infty_{k=1} \Delta_y G(x,  2^k T, y)  (b-u_0) (y) dy.
\eal
\] As shown at the end of Step 1, this series converges uniformly and absolutely. This shows that $u(x, 0)=u_0(x)$ in (\ref{useri}) is equivalent to  that $f$ being given by (\ref{cf}),
completing the proof of the lemma.
\qed

The next lemma shows that if $z=z(x)$ can be reached by the free heat flow at time $T$, then (\ref{areach}) holds with $C_*=\frac{e^+}{T}$ as stated below. The proof follows the idea in \cite{DZ:1} Theorem 2.1.  Since we are dealing with compact manifolds, we are able to reach the explicit constant $C_*$ which is useful for Theorem \ref{thmain}.

\begin{lemma}
\lab{lebhe}
Let $\M$ and $D$ be as in Theorem \ref{thmain} and $u$ a solution of the heat equation
\be
\lab{heq}
\begin{cases} \Delta u -\partial_t u = 0,   \qquad \text{in} \quad D \times (0, T]\\
u(x, t)=0, \quad (x, t) \in \partial D \times (0, T],\\
u(\cdot, 0)=u_0(\cdot) \in L^2(D).
\end{cases}
\ee Then
\be
\lab{dja}
|\Delta^k u(x, T)| \le C \frac{k^{n/4}}{T^{(n+2)/4}} \left(\frac{e}{T}\right)^k k!
\, \Vert u_0 \Vert_{L^2(D)} \le \frac{C}{T^{(n+2)/4}} \left(\frac{e^+}{T}\right)^k k!
\, \Vert u_0 \Vert_{L^2(D)}, \quad k=0, 1, 2, ...
\ee Here
$e$ is $2.71828...$ and $e^+$ is any number greater than $e$.
The constant $C$ depends only on the manifold $\M$ through the dimension, the lower bound of the Ricci curvature, lower bound of the first eigenvalue and volume noncollapsing constant $\inf |B(x_0, 1)|$.
\proof
\end{lemma}

Note
$u^2$, after $0$ extension outside of $D$,  is a subsolution of the heat equation on ${\M} \times (0, T]$. Let $x_0 \in D$ and $k$ be a
positive integer.

 If $T/k \le 1$, then with a suitable translation of time, the
mean value inequality (\ref{pmvi}) with $R=\sqrt{T/k}$,
$\eta=\delta=1/2$ infers that
$$
\al
\sup_{Q_{\sqrt{T/(2k)}} (x_0, T)} u^2 &\le \frac{C_1 (k/T)
}{|B(x_0, \sqrt{T/k})|}
\int_{Q_{\sqrt{T/k}} (x_0, T)} u^2(x, t) \ dx dt\\
&\le \frac{C_2 (k/T)^{(n+2)/2}}{|B(x_0, 1)|}
\int_{Q_{\sqrt{T/k}} (x_0, T)} u^2(x, t) \ dx dt,
\eal
$$ where $Q_r(x_0, T) = B(x_0, r) \times [T-r^2, T]$ is the standard parabolic cube. If the diameter of $\M$ is less than $2 r$, then $B(x, r)$ is regarded as the whole manifold $\M$ here and through out the proof.
In the above we have used the Bishop-Gromov volume comparison
theorem. Note that the above mean value inequality is a local
one since the size of the cubes is less than one. Hence the
constants $C_1$ and $C_2$ are independent of $k$.

If $T/k \ge 1$, then we can apply the mean value inequality on cubes of size $1$ to deduce
\[
u^2(x_0, T) \le \frac{C_2}{|B(x_0, 1)|}
\int_{Q_1(x_0, T)} u^2(x, t) \ dx dt.
\]

Since
$\partial^k_t u$ is also a solution to the heat equation with zero boundary condition, either way it
follows that
\be
\lab{mviqdkt}
\sup_{D \times \{T\}} (\partial^k_t u)^2
\le \frac{C_2 [1+(k/T)^{(n+2)/2}]}{\inf |B(x_0, 1)|}
\int_{D \times [T-T/k, T]} (\partial^k_t u)^2(x, t) \ dx dt.
\ee Next we will bound the right-hand side.

For integers $j=1, 2, \ldots, k,$ consider the space time domains:
$$
\al
\Omega^1_{j}&=  D \times [T (1- \frac{j}{k}), T],\\
\Omega^2_j&= D \times [T (1- \frac{j+0.5}{k}), T].
\eal
$$
Then it is clear that $\Omega^1_j \subset \Omega^2_j \subset
\Omega^1_{j+1}$.

Denote by $\psi^{(1)}_j$ the Lipschitz function of time, which is $0$ on the interval
$[0, -\frac{j+0.5}{k}T]$, $1$ on $[-\frac{j}{k}T, T]$ and linear in between.
Then $|\partial_t \psi^{(1)}_j | \le 2k/T$ a.e.
Since $u$ is a smooth solution to the heat equation, we
 deduce,
 by writing
  $\psi=\psi^{(1)}_j$, that
\[
\al
\int_{\Omega^2_j} &(u_t)^2 \psi \ dx dt = \int_{\Omega^2_j}
u_t \Delta u \psi \ dx dt =-\int_{\Omega^2_j} ((\nabla u)_t \nabla u) \, \psi \ dx dt\\
&= - \frac{1}{2} \int_{\Omega^2_j} (|\nabla u |^2)_t \, \psi
\
dx dt\\
&\le \frac{1}{2} \int_{\Omega^2_j} |\nabla u |^2 \, \psi_t
\
dx dt.
\eal
\]
Therefore,
\be
\lab{omg1j<}
\int_{\Omega^1_j} ( u_t)^2 \ dx dt \le  \frac{k}{T}
\int_{\Omega^2_j} |\nabla u |^2 \ dx dt.
\ee

Denote by $\psi^{(2)}_j$  the Lipschitz function of time, which is $0$ on the interval
$[0, -\frac{j+1}{k}T]$, $1$ on $[-\frac{j+0.5}{k}T, T]$ and linear in between.

 Using $\psi^{(2)}_j u^2$ as a test function in the heat
 equation, the standard Caccioppoli inequality (energy
 estimate)
 between the cubes $\Omega^2_j$ and  $\Omega^1_{j+1}$ shows
 that
 \be
\lab{omg2j<}
\int_{\Omega^2_j} | \nabla u|^2 \ dx dt \le \frac{k}{T}
\int_{\Omega^1_{j+1}}  u^2 \ dx dt.
\ee  A combination of (\ref{omg1j<}) and (\ref{omg2j<}) gives
us
\be
\lab{ddu2r-40}
\int_{\Omega^1_j} ( u_t)^2  \ dx dt \le  (k/T)^2
\int_{\Omega^1_{j+1}} u^2  \ dx dt.
\ee

Since $\partial^j_t u$ is a solution of the heat equation, we can replace $u$ in
(\ref{ddu2r-40}) by $\partial^j_t u$ to deduce, after induction:
\be
\lab{utk<}
\int_{\Omega^1_1} (\partial^k_t u)^2  dx dt \le  (k/T)^{2 k}
\int_{\Omega^1_k} u^2  \ dx dt.
\ee

Note that $\Omega^1_1= D \times [T (1- \frac{1}{k}), T]$ and
$\Omega^1_k=D \times [0, T]$. Substituting
(\ref{utk<}) into (\ref{mviqdkt}), we
find that
\[
 \sup_{D \times \{T\}} (\partial^k_t u)^2
\le \frac{C_2 [1+(k/T)^{(n+2)/2}]}{\inf |B(x_0, 1)|} (k/T)^{2 k}
\int_{D \times [0, T]} u^2(x, t) \ dx dt.
\]
Since $D$ is compact,
\[
\partial_t \int_D u^2(x, t) dx = -2 \int_D |\nabla u|^2(x, t) dx \le - 2 \lambda_1 \int_D u^2(x, t) dx
\]where $\lambda_1>0$ is the first eigenvalue. Therefore
\[
\sup_{D \times \{T\}} (\partial^k_t u)^2
\le \frac{C_2 [1+(k/T)^{(n+2)/2}]}{\inf |B(x_0, 1)| 2 \lambda_1} (k/T)^{2 k}
\int_{D} u^2_0(x) \ dx.
\]Using Stirling's formula
\[
k ! = \sqrt{2 \pi k} \, (k/e)^k \, [1 + O(1/k)],
\]we deduce that
that
\be
\lab{djtu}
|\Delta^k u(x, T)|=|\partial^k_t u(x, T)| \le  C \, \frac{k^{n/4}}{T^{(n+2)/4}} \left( \frac{e}{T} \right)^k \, k! \, \Vert u_0 \Vert_{L^2(D)}
\ee for all integers $k \ge 0$.  Here the constant $C$ depends only on  the manifold $\M$ through the dimension, the lower bound of the Ricci curvature, lower bound of the first eigenvalue and volume noncollapsing constant $\inf |B(x_0, 1)|$.
\qed

Now we are ready to finish
\medskip

\noindent {\it Proof of Theorem \ref{thmain}, part (a).}

Pick $T_0 \in (T-\delta, T]$ with
\be
\lab{deltdef}
\delta=\min \{ \frac{1}{2 A}, \frac{1}{1+ 2 e} T \}
\ee as chosen.
Since, by design, $u(x, T_0)$ is given by the free heat flow on the time interval $[0, T_0]$, i.e.,
\[
u(x, T_0)=\int_D G(x, T_0, y) u_0(y) dy,
\] from Lemma \ref{lebhe}, we know
\[
\al
|\Delta^j u(x, T_0)| &\le C \frac{j^{n/4}}{T_0^{(n+2)/4}} \left(\frac{e}{T_0}\right)^j j! \, \Vert u_0 \Vert_{L^2(D)}\\
&< C_1 \left(\frac{1}{2 T/(1+ 2e)}\right)^j j! \, \Vert u_0 \Vert_{L^2(D)}, \quad j=0, 1, 2, ...
\eal
\] Here we just used the strict inequality
\[
T_0 > T - \frac{1}{1+ 2 e} T = \frac{ 2 e}{1+ 2 e} T.
\]

Recall by assumption
\be
\lab{djzaj}
|\Delta^j z(x)|  \le C A^j j!.
\ee Take
\[
C_* = \max \{ A,  \frac{1}{2 T/(1+ 2e)} \}.
\]Then
\[
|\Delta^j u(x, T_0)| + |\Delta^j z(x)|  \le C_1 C^j_* j!,
\]
\[
T-T_0 < \delta =\frac{1}{2 C_*},
\]which is just condition (\ref{tjie}).
So by applying  Lemma \ref{lemain} with $u(\cdot, T_0)$ as the initial value and $z=z(x)$ as the final value on the time interval $[0, T-T_0]$, we find $v=v(x, t)$ solving
\be
 \lab{cp2}
 \begin{cases}
 \Delta v(x, t) - \partial_t v(x, t)=f(x), \quad (x, t) \in D \times (0, T-T_0], \\
 v(x, t)=0, \quad (x, t) \in \partial D \times (0, T-T_0], \\
 v(x, T-T_0)= z(x), \\
 v(x, 0) = u(x, T_0)=\int_D G(x, t, y) u_0(y) dy,
 \end{cases}
 \ee with
\be
 \lab{cf0}
 f(x) = \int_{\M} \sum^\infty_{k=1} \Delta_y G(x, 2^k (T-T_0), y) \, (u_0-b)(y) dy,
 \ee
 \be
 \lab{fb0}
 b= \sum^\infty_{j=0} \Delta^j z(x) \frac{(T_0-T)^j}{j !}.
 \ee
Take
\[
u(x, t)=
\begin{cases}
v(x, t-T_0), \quad t \in [T_0, T],\\
\int_D G(x, t, y) u_0(y) dy, \quad t \in [0, T_0).
\end{cases}
\]Then $u=u(x, t)$ is the desired solution for
  Theorem \ref{thmain}.

If, in particular,  $z=z(x)$ is reachable by the free heat flow from initial time $0$ to $T$, i.e. (\ref{zgu1}) holds, by Lemma \ref{lebhe}
\[
|\Delta^j z(x)| \le C \frac{j^{n/4}}{T^{(n+2)/4}} \left(\frac{e}{T}\right)^j j! \, \Vert u_0 \Vert_{L^2(D)}, \quad j=0, 1, 2, ...
\] Then we can chose $A=e^+/T$ in the Assumption (\ref{djzaj}). Here $e^+$ is any number greater than $e$.
Therefore the conclusions hold with the following choice of $\delta$ from (\ref{deltdef}):
\[
\delta =\min \{ \frac{1}{2 A}, \frac{1}{1+ 2 e} T \} =\frac{T}{1+2 e},
\] completing the proof of Theorem \ref{thmain}, part (a).
\qed

{\remark   We mention that a condition similar to (\ref{djau0}) on the final state $z=z(x)$ or (\ref{dif0})
on the control function may occur automatically, regardless of the a priori regularity of $u_0$ or $f$. Suppose $u$ is a solution to (\ref{cp}) with final state $u(x, T)=z(x)$ for some stationary control function $f \in L^1(\M)$ and $u_0 \in L^1(\M)$. Then $\partial_t u $ is a solution of the homogeneous heat equation
\[
\begin{cases}
\Delta \partial_t u(x, t) - \partial_t \partial_t u(x, t)=0, \quad (x, t) \in {\M} \times (0, T], \\
\partial_t u(x, 0) = \Delta u_0(x) - f(x).
\end{cases}
\] Here $\Delta u_0$ is understood in the weak sense.  According to Theorem 1.2 in \cite{DZ:1}, $\partial_t u(x, t)$ is analytic  in time for $t>0$ and, there are positive constants $C, C_*$, depending on $T$,  such that
\[
|\Delta^j  \partial_t u(x, t) |_{t=T}| \le C C^j_* j!, \qquad j=0, 1, 2, ....
\] That is
\[
|\Delta^j  (\Delta z(x)-f(x))| \le C C^j_* j!, \qquad j=0, 1, 2, ....
\] In case of null control, i.e., $z(x) \equiv 0$, then
\be
\lab{djfjj}
|\Delta^j  f(x)| \le C C^j_* j!,
\ee showing that (\ref{dif0}) appears automatically.}

It is not hard to make the argument in the remark rigorous and complete

\medskip

\noindent {\it Proof of Theorem \ref{thmain}, part (b).}

We use the method of contradiction. Suppose a solution $u$ to the null control problem exists.
First we prove that  $f$ will be forced to satisfy (\ref{djfjj}).  Since $u$ is a $L^2$ solution, we have
\[
u(x, t) = \int_D G(x, t, y) u_0(y) dy - \int^t_0 \int_{D_0} G(x, t-s, y) f(y) dy ds,
\]which becomes, since $f=f(x)$, that
\be
\lab{u=intd0}
u(x, t) = \int_D G(x, t, y) u_0(y) dy - \int_{D_0} \int^t_0  G(x, s, y) ds f(y) dy.
\ee  Fixing $t>0$ and small $h>0$, by (\ref{u=intd0}), we have
\[
\al
(u(x, t+h)-u(x, t)) h^{-1}& = \int_D ( G(x, t+h, y) - G(x, t, y)) h^{-1} u_0(y) dy \\
&
- \int_{D_0} h^{-1} \int^{t+h}_t G(x, s, y) ds f(y) dy.
\eal
\]By standard heat kernel upper bound we can let $h \to 0$ and apply the dominated convergence theorem to deduce
\[
\partial_t u(x, t) = \int_D \Delta_y  G(x, t, y))  u_0(y) dy - \int_{D_0}G(x, t, y)  f(y) dy.
\] Hence $\partial_t u$ is a smooth solution of the heat equation for $t>0$, with Dirichlet boundary condition.
According to Theorem 1.2 in \cite{DZ:1}, $\partial_t u(x, t)$ is analytic  in time for $t>0$ and
\[
|\Delta^j f(x) | = |\Delta^j  \partial_t u(x, t) |_{t=T}| \le C C^j_* j!, \qquad j=0, 1, 2, ....
\] i.e. (\ref{djfjj}) holds.  This implies $f \in C^\infty_0(D_0)$ and $u$ is smooth for $t>0$.

Fixing $x \in D/D_0$, since $f(x)=0$, we see that $\partial^j_t u(x, T)=0$, for all $j=0, 1, 2, ...$.
But $u(x, t)$ is real analytic for all $t>0$ and continuous down to $t=0$. Hence $u_0(x)=u(x, 0)=0$.
This is a contradiction with the assumption that $u_0$ does not vanish in $D/D_0$, proving Theorem \ref{thmain}, part (b).

\qed

\section{ a control function in a subdomain, positive and negative cases}

In this section we consider the following null control problem:

{\it Let $\omega$ be a nonempty subdomain of $D$ which itself is a bounded, connected domain in a Riemannian manifold or $\textbf{R}^n$. Given $u_0 \in L^2(D)$,
\textbf{find} $g \chi_\omega \in L^2(D \times [0, T])$ and
$u \in L^2([0, T], H^1_0(D)) \cap L^\infty([0, T], L^2(D))$ such that

\be
 \lab{cploc}
 \begin{cases}
 P u(x, t) - \partial_t u(x, t)=g(x, t) \chi_\omega, \quad (x, t) \in D \times (0, T], \\
 u(x, t) = 0, \quad (x, t) \in \partial D \times (0, T],\\
 u(x, T)= 0, \\
 u(x, 0) = u_0(x);\\
 P g(x, t) + \partial_t g(x, t)=0, \quad (x, t) \in D \times (0, T],\\
 g \chi_\omega \in L^2(D \times (0, T]).
 \end{cases}
 \ee  }   Here $P$ is either the Laplace-Beltrami operator $\Delta$ or in case $D$ is a bounded domain in $\textbf{R}^n$,
 \be
 \lab{defopp}
 P= \sum^n_{i, j =1} \partial_i (a_{ij}(x) \partial_j )
 \ee where $\lambda^{-1} I \le (a_{ij}(x)) \le \lambda I$ for some positive constant $\lambda$. In addition, we assume $a_{ij}$ are Lipschitz functions.

 The following is the main result of the section. Since we are not assuming any smoothness of the domain $D$, it is not clear if a solution to problem (\ref{cploc}) exists, due to the lack of observability inequality. See an open problem on p74 of \cite{Zu:1} on null controllability of rough coefficients and domains. Here a progress on rough domains is made by proving that the null control problem (\ref{cploc}) is uniquely solvable on any finite dimensional space generated by eigenfunctions of operator $P$. Moreover, an explicit formula in terms of the eigenvalues and functions is given. See also Remark \ref{bmcoef} for a negative answer  when $(a_{ij})$ is just bounded measurable.

 \begin{theorem}
\lab{thloccon}
 Let $\M$ be a n dimensional, compact Riemannian manifold without boundary or $ \textbf{R}^n$. Let $D \subset \M$ be a bounded,  connected domain. Let $G=G(x, t, y)$ be the Dirichlet heat kernel of the operator $P$ on domain $D$.

(a). The null control problem is equivalent to the following integral equation for the unknown function $\phi$ in the standard Hilbert space $\textbf{H}$ in (\ref{deHspace}):
\be
\lab{inteqphi}
\int_D \left[ \int^T_0 \int_\omega G(x, s, z) G(z, s, y) dz ds \right] \, \phi(y) dy = \psi(x)
\equiv \int_D G(x, T, y) u_0(y) dy.
\ee  If $\phi$ is found, then the control function is given by
$
g(x, s)=\int_D G(x, T-s, y) \phi(y) dy.
$

(b).  Let $\lambda_j<0$, $j=1, 2, ...$,  be all the eigenvalues of $P$ on $D$ counting multiplicity and $\eta_j$ be the $L^2$ normalized and orthonormal eigenfunctions. Define
\be
\al
\alpha_{ij} &= \frac{1-e^{(\lambda_i+\lambda_j) T}}{|\lambda_i + \lambda_j|} <\eta_i, \eta_j>_{L^2(\omega)}, \quad
 \,  \, \beta_i =\int_D \psi(y) \eta_i(y) dy=e^{\lambda_i T} <\eta_i, u_0>_{L^2(D)}\\
 s_j &= \int_D \phi(y)  \eta_j(y) dy.
\eal
\ee The integral equation (\ref{inteqphi}) is equivalent to the infinite dimensional linear algebraic system for the unknown $s_j$:
\be
\sum^\infty_{j=1} \alpha_{ij} \,  s_j = \beta_i.
\ee

(c). For any positive integer $m$, the matrix $A_m \equiv (\alpha_{ij})^{m}_{i, j=1}$ is positive definite and hence invertible. Let $V_m$ be the linear space spanned by $\{ \eta_j \}^m_{j=1}$.

 Given any initial value $u_0 \in V_m$, the null control problem (\ref{cploc}) has a unique solution given by the explicit formula: $g(x, s)=\int_D G(x, T-s, y) \phi(y) dy =\sum^m_{j=1} e^{\lambda_j (T-s)} s_j \eta_j(x)$,
 \[
 \phi = \sum^m_{j=1} s_j \eta_j(x)
 \]with $s_j = \sum^m_{i=1} \theta_{ij} \beta_i$ where $(\theta_{ij})$ is the inverse matrix of $A_m$.
\proof
\end{theorem}

  The proof is divided into 4 steps.
\medskip

{\it Step 1.} set up a pertinent Hilbert space and an operator.

 Recall $G=G(x, t, y)$ be the Dirichlet heat kernel of the operator $P$ on domain $D$.
 Given $\phi \in L^2(D)$, let
 \[
 h=h(x, t) =\int_D G(x, t, y) \phi (y) dy \equiv e^{t P} \phi
 \]i.e. $h$ is the solution of the initial boundary value problem of the heat type equation:
 \be
 \lab{hehphi}
 \begin{cases}
 P h(x, t) - \partial_t h(x, t)=0, \quad (x, t) \in D \times (0, T], \\
 h(x, t) = 0, \quad (x, t) \in \partial D \times (0, T],\\
 h(x, 0) = \phi(x).
 \end{cases}
 \ee
 Consider the Hilbert space $\textbf{H}$ which is the completion of the set
\be
\lab{deHspace}
 \{ \phi \in L^2(D)  \, | \,  \int^T_0 \int_\omega \left| e^{s P} \phi (y)\right|^2 dy ds <\infty \}
\ee under the norm
\be
\lab{defnormH}
\Vert \phi \Vert_\textbf{H} =  \left( \int^T_0 \int_\omega \left| e^{s P} \phi (y)\right|^2 dy ds )\right)^{1/2}  = \Vert h \Vert_{L^2(\omega \times [0, T])}.
\ee By the uniqueness result is \cite{Lin:1} p136, if $h=0$ on $\omega \times [0, T]$ then
$h=0$ on the whole domain $D \times [0, T]$.  This is the only place where we used the assumption that
the leading coefficients $a_{ij}$ are Lipschitz.

We will only give a proof for the case $P=\Delta$, the Laplace-Beltrami operator.  For the other case, i.e. $P$ is a uniformly elliptic operator with Lipschitz coefficients, we just need to replace the heat kernel upper bound on a compact manifold by the bound for $G(x, t, y)$ in Aronson's classic paper \cite{Ar:1}.  No smoothness of the heat kernel is needed since solutions of (\ref{cploc}) and (\ref{hehphi}) are understood in $L^2([0, T], H^1_0(D)) \cap L^\infty([0, T], L^2(D)) $ sense.

The main idea is to consider a kernel function arising from Duhamel's formula and use it to define a bilinear form on the Hilbert space $ \textbf{H}$.

Define a kernel function
\be
\lab{defkk}
K(x, T, y) \equiv \int^T_0 \int_\omega G(x, s, z) G(z, s, y) dz ds.
\ee From standard heat kernel upper bound
\be
\lab{ggub}
G(x, t, y) \le C_1\left( \frac{1}{t^{n/2}} +1 \right) e^{- C_2 d^2(x, y)/t}
\ee and routine integration on the bounded domain $\omega$, we know that
\be
\lab{pbkk}
0 \le K(x, T, y) \le \frac{C_T}{d(x, y)^{n-2}}.
\ee Here $C_T$ is a positive constant depending on $T$, $\M$;  $d(x, y)$ is the geodesic distance for $x, y \in \M$.
 It is also clear that $K$ is symmetric in $x, y$ due to the similar symmetry of the heat kernel.
Note that if $\omega$ is the full domain $D$, then the reproducing property of the heat kernel shows
that $K(x, T, y) = \int^T_0  G(x, 2 s, y) ds$.  Anyway the singularity of kernel $K$ is dominated by that of the Green's function of the Laplacian or operator $P$. Since $D$ is bounded, $K$ is an integrable kernel, so the operator
\be
\lab{defopl}
L  \phi (x)  \equiv \int_D K(x, T, y) \phi (y) dy
\ee maps $L^2(D)$ into itself.  Moreover,  for any $\phi, \psi \in L^2(D)$, we have, after using Fubini theorem and Cauchy-Schwarz inequality, that
\be
\lab{lphisi}
\al
< L\phi, \psi >_{L^2} &\equiv \int_D L\phi(x)  \, \psi(x) dx\\
&=\int^T_0 \int_D \int_D  \int_\omega G(x, s, z) G(z, s, y) dz \, \phi(y) dy \, \psi(x) dx  ds\\
&=\int^T_0 \int_\omega \int_D G(z, s, x)  \psi(x) dx \,  \int_D  G(z, s, y) \phi(y) dy \,  dz  ds\\
&\equiv < \psi, \phi>_{\textbf{H}}\\
&\le \left( \int^T_0 \int_\omega | e^{s \Delta} \phi(z)|^2 dz ds \right)^{1/2}
\left( \int^T_0 \int_\omega | e^{s \Delta} \psi(z)|^2 dz ds \right)^{1/2}
\eal
\ee This implies that
\be
\lab{lyouji}
|< L\phi, \psi >_{L^2}| = | < \psi, \phi>_{\textbf{H}} |  \, \le \Vert \phi \Vert_{\textbf{H}} \, \Vert \psi \Vert_{\textbf{H}}.
\ee Also, in case $\phi=\psi$, the 3rd line of (\ref{lphisi}) says that
\be
\lab{lqz}
< L\phi, \phi >_{L^2} = \Vert \phi \Vert^2_{\textbf{H}}.
\ee

Given $\phi \in L^2(D)$, by standard energy estimate for the heat equation, there is a positive constant $C_T>0$ such that
\[
\al
\Vert L \phi \Vert^2_{\textbf{H}} &\le C_T \Vert L \phi \Vert^2_{L^2} = C_T <L \phi, \, L \phi>_{L^2}\\
&=C_T <\phi, \,  L \phi>_{\textbf{H}}, \qquad \text{by} \quad (\ref{lyouji}),\\
& \le C_T \Vert  \phi \Vert_{\textbf{H}} \, \Vert L \phi \Vert_{\textbf{H}}
\eal
\]This implies
\be
\lab{lhjie}
\Vert L \phi \Vert_{\textbf{H}} \le C_T \Vert  \phi \Vert_{\textbf{H}},
\ee i.e. $L$ can be extended to a bounded linear operator from $\textbf{H}$ to itself.
From now on, we will just use $L$ to denote this expanded operator. For example, when we write $L \phi$, it makes sense for all $\phi \in \textbf{H}$ which may be much larger than $L^2(D)$.

\medskip

{\it Step 2.} Converting (\ref{cploc}) into an equation involving operator $L$.

Now suppose for some $\psi \in \textbf{H}$, there exists $\phi \in \textbf{H}$ such that
\be
\lab{lphi=si}
L \phi = \psi.
\ee From the definition of $L$ in (\ref{defopl}) and (\ref{defkk}), this identity is equivalent to:
\be
\lab{ggphis}
\int_D \int^T_0 \int_\omega G(x, s, z) G(z, s, y) dz ds \, \phi(y) dy = \psi(x).
\ee As explained at the end of Step 1, the above identity is understood in the $\textbf{H}$ space sense.

Making a change of time variable $s \to T -s$, (\ref{ggphis}) becomes
\be
\lab{ggphis2}
 \int^T_0 \int_\omega G(x, T-s, z) \int_D  G(z, T-s, y)  \phi(y) dy dz ds = \psi(x).
\ee

 Next we take
\be
\lab{psi=z-}
\psi(x) = \int_D G(x, T, y) u_0(y) dy.
\ee Then $\psi \in \textbf{H}$ since it is in $L^2(D)$ and the solution of the heat equation with $\psi$ as initial value is $L^2$ in $D$.
Now suppose, with this choice of $\psi$, the equation (\ref{lphi=si}) has a solution $\phi \in \textbf{H}$. Then
take
\be
\lab{defgphi}
g(z, s) \equiv  \int_D  G(z, T-s, y)  \phi(y) dy
\ee with $g \chi_\omega \in L^2(\omega \times [0. T])$. The meaning of (\ref{defgphi}) needs some explanation. We have no a priori knowledge if $\phi$ is integrable on $D$ or not. What it means is that there exists a sequence $\phi_j \in L^2(D)$ such that $\Vert \phi_j - \phi \Vert_{\textbf{H}} \to 0$ i.e.
$\int_D  G(z, T-s, y)  \phi_j(y) dy$ converges to $g$ is $L^2(\omega \times [0, T])$.

Substitution of (\ref{defgphi}) and (\ref{psi=z-}) into (\ref{ggphis2}) gives us
\be
0 = \int_D G(x, T, y) u_0(y) dy -  \int^T_0 \int_\omega G(x, T-s, z) g(z, s) \chi_w  dz ds,
\ee which shows that the function
\be
\lab{uxtlc}
u(x, t) \equiv  \int_D G(x, t, y) u_0(y) dy -  \int^t_0 \int_\omega G(x, t-s, z) g(z, s) \chi_w  dz ds
\ee solves (\ref{cploc}). Indeed, since $g \chi_\omega \in L^2([0, T] \times D)$, the standard heat kernel or energy estimate on
(\ref{uxtlc}) says $u \in L^2([0, T], H^1_0(D)) \cap L^\infty([0, T], L^2(D))$.
This  proves part (a) of the theorem, i.e. the equivalence of (\ref{lphi=si}) with problem (\ref{cploc}).  In fact it is also known that
$u \in C([0, T], L^2(D))$. See Remark \ref{rmkuc} below.
\medskip

{\it Step 3. conversion to an infinite dimensional linear system}
\medskip

First let us remark that $L$ is injective. There are different quick proof of this fact. For example we take $u_0=0$ and suppose $g$ is a solution to problem (\ref{cploc}).  From Step 2 with $\psi =0$ in (\ref{lphi=si}), if there exists $\phi \in \textbf{H}$ such that
\[
L \phi =0.
\]Then $\forall \eta_j \in \textbf{H} \cap L^2(D)$, we deduce, from (\ref{lyouji}) that
\[
0= < L \phi,  \eta_j >_{L^2} = <\phi, \eta_j>_{\textbf{H}}.
\] Selecting $\eta_j$ such that $\eta_j \to \phi$ in $\textbf{H}$ norm, we conclude
\[
\Vert \phi \Vert_{\textbf{H}}=0.
\]This shows that $L$ is injective

  From Step 2 with $\psi =\int G(x, T, y) u_0(y) dy$ in (\ref{lphi=si}), we next give a formal solution to
\be
\lab{lphi=si2}
L \phi =\psi
\ee in the form of a infinite dimensional linear system.

Let $\lambda_j<0$, $j=1, 2, ...$,  be all the eigenvalues of $\Delta$ on $D$ counting multiplicity and $\eta_j$ be the $L^2(D)$ normalized and orthonormal eigenfunctions. It is well known that the heat kernel is given by
\be
\lab{hkxilie}
G(x, T, y) = \sum^\infty_{j=1} e^{\lambda_j T} \eta_j(x) \eta_j(y).
\ee Note the series on the right hand side converges pointwise absolutely for each fixed $T$ due to
\[
\al
&\left| \sum^\infty_{j=1} e^{\lambda_j T} \eta_j(x) \eta_j(y) \right|
\le \left( \sum^\infty_{j=1} e^{\lambda_j T} \eta^2_j(x) \right)^{1/2}
\left( \sum^\infty_{j=1} e^{\lambda_j T} \eta^2_j(y) \right)^{1/2}\\
&=\sqrt{G(x, T, x) G(x, T, y)} \le \frac{C_T}{T^{n/2}}.
\eal
\]
Hence, after using dominated convergence theorem and heat kernel upper bound (\ref{ggub}),
\[
G(x, s, z) G(z, s, y) \le C_1\left( \frac{1}{s^{n/2}} +1 \right)^2 e^{- C_2 d^2(x, z)/s - C_2 d^2(z, y)/s},
\] which is integrable on $D \times [0, T]$,  we find that the kernel $K$ of the operator $L$ is given by
\be
\lab{kerxilie}
\al
K(x, T, y)&=\int^T_0 \int_\omega G(x, s, z) G(z, s, y) dz ds\\
&=\int^T_0 \int_\omega \lim_{m \to 0} \left( \sum^m _{j=1} e^{\lambda_j s} \eta_j(x) \eta_j(z) \,
\sum^m _{j=1} e^{\lambda_j s} \eta_j(z) \eta_j(y) \right) dzds\\
& = \sum^\infty_{i, j=1} \frac{1-e^{(\lambda_i+\lambda_j) T}}{|\lambda_i + \lambda_j|} <\eta_i, \eta_j>_{L^2(\omega)} \eta_i(x) \eta_j(y).
\eal
\ee Then equation (\ref{lphi=si2}) can be transformed into the infinite dimensional linear system:
\be
\sum^\infty_{i, j=1} \frac{1-e^{(\lambda_i+\lambda_j) T}}{|\lambda_i + \lambda_j|} <\eta_i, \eta_j>_{L^2(\omega)} \int_D \phi(y)  \eta_j(y) dy  \, \eta_i(x) = \sum^\infty_{i=1} \int_D \psi(y) \eta_i(y) dy \, \eta_i(x).
\ee This is equivalent to
\be
\lab{inftyfcz}
\sum^\infty_{i, j=1} \frac{1-e^{(\lambda_i+\lambda_j) T}}{|\lambda_i + \lambda_j|} <\eta_i, \eta_j>_{L^2(\omega)} \int_D \phi(y)  \eta_j(y) dy = \int_D \psi(y) \eta_i(y) dy.
\ee  Write
\be
\alpha_{ij} = \frac{1-e^{(\lambda_i+\lambda_j) T}}{|\lambda_i + \lambda_j|} <\eta_i, \eta_j>_{L^2(\omega)}, \quad
s_j = \int_D \phi(y)  \eta_j(y) dy,  \, \beta_i =\int_D \psi(y) \eta_i(y) dy.
\ee Then we can write (\ref{inftyfcz}) in the compressed form
\be
\sum^\infty_{j=1} \alpha_{ij} \,  s_j = \beta_i.
\ee This is equivalent to (\ref{lphi=si2}), proving part (b) of the theorem.

{\it Step 4. null controllability on any finite dimensional space spanned by eigenfunctions}
\medskip

Let $m$ be any positive integer. We claim that the matrix $A_m \equiv (\alpha_{ij})^{m}_{i, j=1}$ is positive definite and hence invertible. Here goes the proof. Define an approximate heat kernel to be
\be
\lab{hkxiliem}
G_m(x, T, y) = \sum^m_{j=1} e^{\lambda_j T} \eta_j(x) \eta_j(y).
\ee Similarly define an approximate $K$ kernel to be
\be
\lab{kxiliem}
\al
K_m(x, T, y)&=\int^T_0 \int_\omega G_m(x, s, z) G_m(z, s, y) dz ds\\
& = \sum^m_{i, j=1} \frac{1-e^{(\lambda_i+\lambda_j) T}}{|\lambda_i + \lambda_j|} <\eta_i, \eta_j>_{L^2(\omega)} \eta_i(x) \eta_j(y)\\
&=\sum^m_{i, j=1} \alpha_{ij} \, \eta_i(x) \eta_j(y).
\eal
\ee Let $V_m$ be the linear space spanned by $\{\eta_j \}^m_{j=1}$. Consider the operator $L_m:
V_m \to V_m$:
\be
\lab{defoplm}
L_m \phi = \int_D K_m(x, T, y) \phi (y) dy.
\ee Then, similar to (\ref{lyouji}), we have, for $\phi= \sum^m_{j=1} b_j \, \eta_j(x)$,
\be
\lab{aij=int}
\sum^m_{i, j=1} \alpha_{ij} b_i b_j = <L_m \phi, \phi>_{L^2(D)} =\int^T_0 \int_\omega \left(
\int_D G_m(x, s, z) \phi(z) dz \right)^2 dxds.
\ee Therefore $A_m$ is a semi-positive definite matrix.  Since $L_m$ is the restriction of the operator $L$ on the space $V_m$ and $L$ is injective from the beginning of Step 3,  we have proved the claim that $A_m$ is positive definite and invertible.
Hence the system
\be
\lab{as=bm}
\sum^m_{j=1} \alpha_{ij} \,  s_j = \beta_i.
\ee has a unique solution
\be
s_j =\sum^m_{i=1} \theta_{ij}  \beta_i, \qquad j=1, ..., m,
\ee where $(\theta_{ij})$ is the inverse matrix of $A_m=(\alpha_{ij})^m_{i, j=1}$.

Given any $u_0 \in V_m$, we take
\be
\lab{psibtsj}
\al
\psi&=\int_D G_m(x, T, y) u_0(y) dy, \quad \beta_i = <\psi, \eta_i>_{L^2(D)}=e^{\lambda_i T} <u_0, \eta_i>_{L^2(D)}, \quad\\
&\phi = \sum^m_{j=1} s_j \eta_j.
\eal
\ee Then (\ref{as=bm}) implies
\[
L_m \phi = \psi
\] which is equivalent to
\[
0 = \int_D G_m(x, T, y) u_0(y) dy -  \int^T_0 \int_\omega G_m(x, T-s, z) g_m(z, s) \chi_w  dz ds
\]where
\be
\lab{defgm}
g_m(z, s) =\int_D G_m(z, T-s, y) \phi(y) dy.
\ee Hence
\[
u=\int_D G_m(x, t, y) u_0(y) dy -  \int^t_0 \int_\omega G_m(x, t-s, z) g_m(z, s) \chi_w  dz ds
\]is a solution to the null control problem (\ref{cploc}).

This  completes the proof of the theorem.
\qed

{\remark
\lab{bmcoef}

On p74 \cite{Zu:1}, a question is raised on the null controllability of the equation in (\ref{cploc}) when the leading coefficients $(a_{ij})$ are just bounded measurable. Here the only requirement for the control function is that it is $L^2(\omega \times [0, T])$.  It turns out that the answer is {\it negative} even for some $(a_{ij})$ which are H\"older continuous, due to the example in the paper \cite{Fi:1}. There the author constructed an elliptic matrix $(a_{ij})$ which are H\"older continuous such that the equation
\[
\sum^n_{i, j=1} \partial_i (a_{ij}(x) \partial_j u(x) ) = \lambda u(x)
\] has a nontrivial solution $u \in C^\infty_0(\textbf{R}^n)$. Then $v \equiv e^{\lambda t} u(x)$ is a solution of the corresponding parabolic equation for which the observability inequality
\be
 \int_D v^2(x, T) dx \le C_T \int^T_0 \int_\omega v^2(x, t) dxdt
\ee obviously fails when $\omega$ is disjoint from the
support of $u$.  Since the validity of the observability inequality for all finite $L^2(D)$ solutions of the heat equation is equivalent to null controllability (see Section 4.3 \cite{Zu:1} e.g.), the later also fails for some initial values.

We wish to thank Prof. F. H. Lin for informing me the paper \cite{Fi:1} and for the above construction of a solution of the parabolic equation.
}

{\remark
For some special domains such as a box in $\textbf{R}^n$, the eigenvalues and functions are explicitly known and the formula in the theorem becomes easily computable by computer when $P=\Delta$ in $\textbf{R}^n$.

There is a chance to prove null controllability in general if $\beta_i$ decays sufficiently fast.
Indeed from (\ref{psibtsj}) and (\ref{defgm}), we have
\[
g_m(x, t) = \sum^m_{j=1} e^{\lambda_j (T-t)} s_j \eta_j(x).
\]Hence
\[
\al
\Vert g_m \Vert^2_{\textbf{H}}& = \int^T_0 \int_\omega g^2_m(x, t) dxdt \\
&=\sum^m_{i, j=1}  \frac{1-e^{(\lambda_i+\lambda_j) T}}{|\lambda_i + \lambda_j|} <\eta_i, \eta_j>_{L^2(\omega)} s_i s_j\\
&=\sum^m_{i, j=1} \alpha_{ij} s_i s_j = \sum^m_{j=1} \beta_j s_j.
\eal
\] Note that since $u_0 \in L^2(D)$, one can prove that $\beta_j$ decays like $e^{-|\lambda_j| T^-}$, where
$T^-$ is any number strictly less than $T$. So, one would have a uniform bound for $\Vert g_m \Vert_{\textbf{H}}$ if one can prove that $s_j$ grows slower than $ e^{|\lambda_j| T^-}$. Then weak sub-convergence of $g_m$ in $\textbf{H}$ follows.
}

{\remark Using $g_m$ as a test function for the equation of $u$ in (\ref{cploc}), we obtain
\[
\int^T_0 \int_\omega g^2_m(x, t) dxdt = \int_D u_0(x) g_m(x, 0) dx.
\] In case the observability inequality holds:
\[
\int_D g^2_m(x, 0) dx \le C_T \int^T_0 \int_\omega g^2_m(x, t) dxdt
\] as in the case $D$ is smooth, we can obtain the uniform bound
\[
\int^T_0 \int_\omega g^2_m(x, t) dxdt \le C C_T \int_D u^2_0(x) dx.
\]This implies weak sub-convergence of $\{ g_m \}$ in $\textbf{H}$ and provides an effective computation method for the control function without finding the minimizer of the traditional $J$ functional.
}

{\remark
\lab{rmkuc}

The following known fact is proven for completeness, which is useful in justifying that the final state $u(x, T)=0$.

 If $u$ is a solution to (\ref{cploc}), then $u \in C([0, T], L^2(D))$ i.e. it is continuous in time in the  $L^2(D)$ sense for the following reason.   Given $t_1, t_2 \in (0, T]$, from (\ref{uxtlc}), after a shift in time, we deduce
\[
\al
u(x, t_2) - u(x, t_1)& = \int_D G(x, t, y)\big{|}^{t_2}_{t_1} u_0(y) dy
- \int^{t_2}_{t_1} \int_D G(x, s, y) g(y, t_2-s) \chi_\omega dyds\\
&\qquad  - \int^{t_1}_0 \int_D G(x, s, y) \left(g(y, t_2-s) -g(y, t_1-s) \right) \chi_\omega dyds.
\eal
\]
From here, for any smooth test function $\eta$, applying
 Young's inequality, we find
\[
\al
&\left| \int_D (u(x, t_2) - u(x, t_1)) \eta(x) dx \right| \\
&\le \sup_{y \in D} \Vert G(\cdot, t_2, y) -G(\cdot, t_1, y) \Vert_{L^1(D)} \Vert u_0 \Vert_{L^2(D)}  \Vert \eta \Vert_{L^2(D)} + \sqrt{t_2-t_1} \Vert g \Vert_{L^2([t_1, t_2] \times \omega)} \, \Vert \eta \Vert_{L^2(D)}\\
&\qquad +
 \sqrt{t_1} \Vert g(\cdot, t_2-\cdot) - g(\cdot, t_1-\cdot) \Vert_{L^2([0, t_1] \times \omega)} \, \Vert \eta \Vert_{L^2(D)}\\
&=o(1) \to 0, \qquad \text{as} \quad t_1 \to t_2.
\eal
\] In the above we have used the fact that $\int_D G(x, t, y) dy \le 1$. This yields the $L^2$ continuity of $u$ in $t$.
}

\section{ An inversion formula for the heat kernel}

In this short section, we present a byproduct of the proof of the theorem, an inversion formula for the heat kernel. Due to the ubiquity of the heat kernel, this formula may be of independent interest.

\begin{proposition}
\lab{prinv}
 Let $\M$ be a n dimensional, compact Riemannian manifold without boundary and $D \subset \M$ a smooth domain or $D=\M$.
 Suppose $u=u(x, t)$ is a solution of following problem with $ u_0 \in L^2(D)$ as the initial value:
 \be
 \lab{ibvp}
 \begin{cases}
 \Delta u(x, t) - \partial_t u(x, t)=0, \quad (x, t) \in D \times (0, T], \\
 u(x, t)=0, \quad (x, t) \in \partial D \times (0, T], \\
 u(x, 0) = u_0(x).
 \end{cases}
 \ee Then for all $t \in ( (1-e^{-1})T, T]$, we have
 \be
 \lab{hkinv}
 u(x, t) =  \left( \sum^\infty_{k=0} \frac{(t-T)^k}{k !} \Delta^k \right) u(x, T).
 \ee
\end{proposition}

{\remark
Note that
\[
u(x, T) = \int_D G(x, T; y, t) u(y, t) dy \equiv P_{t, T} (u(\cdot, t)) = e^{(T-t) \Delta} u(\cdot, t)
\]where $G=G(x, T; y, t)=G(x, T-t, y)$ is the heat kernel with $0$ boundary value.
The heat operator $P_{t, T}$ can be regarded as a linear operator from $L^2(D)$ to $C^\infty(D) \cap C_0(D)$.
The Proposition shows that the operator
\be
\lab{hopinv}
\sum^\infty_{j=0} \frac{(t-T)^j}{j !} \Delta^j
\ee is the inverse of $P_{t, T}$ in a subspace of $C^\infty(D) \cap C_0(D)$. Note (\ref{hopinv}) is nothing but $e^{(t-T) \Delta}$ formally.  The point of the proposition is that it is well defined on the space of functions reachable by the free heat flow in a fixed time interval.

In case $u$ is also analytic in space variables, as in the Euclidean setting, this formula provides a way to recover a state in the past from information at one point for the state at present. This kind of formula has been sought after in the inverse problem community.}
\medskip

{\noindent \it Proof of Proposition \ref{prinv}}

According to Lemma \ref{lebhe}, the following inequality holds

\[
|\Delta^k u(x, T)| \le C \frac{k^{n/4}}{T^{(n+2)/4}} \left(\frac{e}{T}\right)^k k!
\, \Vert u_0 \Vert_{L^2(D)}, \quad k=0, 1, 2, ...
\] If $t \in ( (1-e^{-1})T, T]$, then
\[
0 \le T-t < T/e.
\]Hence the series on the right hand side of (\ref{hkinv}):
\be
\lab{hkinvr}
  \sum^\infty_{k=0} \frac{(t-T)^k}{k !} \Delta^k  u(x, T),
\ee and its formal derivatives with respect to $t$
and
\[
\sum^\infty_{k=0} \frac{(t-T)^k}{k !} \Delta^k  \Delta u(x, T)
\] all converge uniformly and absolutely. Then it is clear that $(\ref{hkinvr})$ is a solution of the heat equation with final value $u(x, T)$. The uniqueness of the backward heat equation shows that
$(\ref{hkinv})$ is true.   \qed

{\remark  The result can be extended to some noncompact setting with a different time intervals. By iterating formula (\ref{hkinv}), the  time to which the heat kernel can be inverted may be improved to any $t \in (0, T]$.}

\medskip
{\bf Acknowledgment.} We wish to thank Professor H. J. Dong for suggesting a shorter way to get (\ref{cf}) and finding an error in an earlier version,
Professor F. H. Lin for adding two important references, and Professor X. H. Pan for making corrections.

\bigskip

\noindent e-mails:
qizhang@math.ucr.edu


\begin{thebibliography}{00}

\bibitem[Ar]{Ar:1} Aronson, D. G. {\it Non-negative solutions of linear parabolic equations.} Ann. Scuola Norm. Sup. Pisa Cl. Sci. (3) 22 (1968), 607-694.


\bibitem[ABGM]{ABGM:1} Ammar Khodja, Farid; Benabdallah, Assia; Gonzalez-Burgos, Manuel; Morancey, Morgan,{\it
Quantitative Fattorini-Hautus test and minimal null control time for parabolic problems.}
J. Math. Pures Appl. (9) 122 (2019), 198-234.


\bibitem[CJ]{CJ:1}
Christensen, Ann-Eva;  Johnsen, Jon, {\it  Final value problems for parabolic differential equations and their well-posedness}, arXiv:1707.02136, Axioms, 7(2), 1-36, 2018.

\bibitem [DM]{DM:1} Duyckaerts, Thomas; Miller, Luc, {\it
Resolvent conditions for the control of parabolic equations.}
J. Funct. Anal. 263 (2012), no. 11, 3641-3673.

\bibitem [DP]{DP:1}  H. J.  Dong and X. H. Pan, {\it
Time analyticity for inhomogeneous parabolic equations and the Navier-Stokes equations in the half space},
arXiv:2003.12915.

 \bibitem [DZ]{DZ:1}  H. J.  Dong and Q. S. Zhang, {\it Time analyticity for the heat equation and Navier-Stokes equations} ,arXiv:1907.01687, J. Funct. Anal. https://doi.org/10.1016/j.jfa.2020.108563, 2020.

\bibitem[EMZ]{EMZ:1} L. Escauriaza, S. Montaner and C. Zhang, {\it Analyticity of solutions to parabolic evolutions and applications}, SIAM J. Math. Anal. 49 (2017), no. 5, 4064-4092.

\bibitem [EZ]{EZ:1} Ervedoza, Sylvain; Zuazua, Enrique, {\it Sharp observability estimates for heat equations.} Arch. Ration. Mech. Anal. 202 (2011), no. 3, 975-1017.


\bibitem [Fi]{Fi:1} Filonov, N. {\it  Second-Order Elliptic Equation of Divergence Form Having a Compactly Supported Solution}, Journal of Mathematical Sciences volume 106, p3078-3086 (2001)

\bibitem [FI]{FI:1}  Fursikov, A. V.; Imanuvilov, O. Yu. {\it Controllability of evolution equations. } Lecture Notes Series, 34. Seoul National University, Research Institute of Mathematics, Global Analysis Research Center, Seoul, 1996. iv+163 pp.



\bibitem [FR]{FR:1}Fattorini, H. O.; Russell, D. L. {\it
Exact controllability theorems for linear parabolic equations in one space dimension.}
Arch. Rational Mech. Anal. 43 (1971), 272-292.


\bibitem [Li] {Li:1} Li, Peter,
{\it Geometric analysis.} Cambridge Studies in Advanced Mathematics, 134. Cambridge University Press, Cambridge, 2012. x+406 pp.

\bibitem [Lin] {Lin:1} Lin, Fang-Hua, {\it A uniqueness theorem for parabolic equations.} Comm. Pure Appl. Math. 43 (1990), no. 1, 127-136.

\bibitem [Lin2] {Lin:2} Lin, Fanghua, {\it  Remarks on a backward parabolic problem.} Methods Appl. Anal. 10 (2003), no. 2, 245-252.

\bibitem [Lio] {Lio:1} Lions, J.-L. {\it Optimal control of systems governed by partial differential equations.} Translated from the French by S. K. Mitter. Die Grundlehren der mathematischen Wissenschaften, Band 170 Springer-Verlag, New York-Berlin 1971.

\bibitem [Lp] {Lp:1} Lions, P.-L. {\it Optimal control of diffusion processes and Hamilton-Jacobi-Bellman equations.} I. The dynamic programming principle and applications. Comm. Partial Differential Equations 8 (1983), no. 10, 1101-1174


\bibitem [LL] {LL:1} Le Rousseau, J\'er\^ome; Lebeau, Gilles, {\it On Carleman estimates for elliptic and parabolic operators.} Applications to unique continuation and control of parabolic equations. ESAIM Control Optim. Calc. Var. 18 (2012), no. 3, 712-747.



\bibitem [LR] {LR:1} Lebeau, G.; Robbiano, L. {\it Contr\^ole exact de l'\'equation de la chaleur. (French) [Exact control of the heat equation]} Comm. Partial Differential Equations 20 (1995), no. 1-2, 335-356.
\bibitem [LT]{LT:1}  Lasiecka, Irena; Triggiani, Roberto, {\it Control theory for partial differential equations}, I, II,  Cambridge University Press, Cambridge, 2000.

\bibitem [LZ]{LZ:1}   Lebeau, Gilles; Zuazua, Enrique, {\it Null-controllability of a system of linear thermoelasticity}. Arch. Rational Mech. Anal. 141 (1998), no. 4, 297-329.


\bibitem [LZZ] {LZZ:1} Lopez, Antonio; Zhang, Xu; Zuazua, Enrique, {\it
Null controllability of the heat equation as singular limit of the exact controllability of dissipative wave equations.}
J. Math. Pures Appl. (9) 79 (2000), no. 8, 741-808.

\bibitem [MZ] {MZ:1} Micu, Sorin; Zuazua, Enrique, {\it Regularity issues for the null-controllability of the linear 1-d heat equation.} Systems Control Lett. 60 (2011), no. 6, 406-413.



\bibitem [Tr]{Tr:1} Tr\"oltzsch, Fredi, {\it Optimal control of partial differential equations. Theory, methods and applications}.  Graduate Studies in Mathematics, 112. American Mathematical Society, Providence, RI, 2010. xvi+399 pp

\bibitem [Ru]{Ru:1}  Russell, David L. {\it Controllability and stabilizability theory for linear partial differential equations: recent progress and open questions.} SIAM Rev. 20 (1978), no. 4, 639-739.


\bibitem [RW]{RW:1}  Russell, David L.; Weiss, George,
{\it A general necessary condition for exact observability.}
SIAM J. Control Optim. 32 (1994), no. 1, 1-23.



\bibitem [Z]{Z:1}  Q. S. Zhang, {\it A note on time analyticity for ancient solutions  of the heat equation}, Proc.Amer. Math. Soc. 148 (2020), no. 4, 1665-1670.

\bibitem [Zu]{Zu:1} Zuazua, Enrike, {\it Controllability and observability of partial differential equations: some results and open problems.} Handbook of differential equations: evolutionary equations. Vol. III, 527-621, Handb. Differ. Equ., Elsevier/North-Holland, Amsterdam, 2007.

\end{thebibliography}
\enddocument